\documentclass[12pt]{article}
\usepackage[utf8]{inputenc}
\usepackage{amssymb,amsmath}
\usepackage{enumerate}
\usepackage[amsmath,amsthm,thmmarks]{ntheorem}
\usepackage{ulem}
\usepackage[utf8]{luainputenc}
\usepackage{geometry}
\usepackage[pdftex]{graphicx}
\usepackage{amstext}
\usepackage{amsmath}
\usepackage[T1,T2A]{fontenc}
\usepackage[utf8]{inputenc}
\usepackage[english,main=russian]{babel} 
\usepackage{mathtools}
\usepackage{colordvi}
\usepackage{amsfonts}

\def\F{{\cal F}}
\def\repV{\Upsilon_{\rm v}}
\def\repVSym{\Upsilon_{\rm v}^{\rm sym}}

  \newtheorem{theorem}{Теорема}[section]
  \newtheorem{lemma}{Лемма}[section]
  \newtheorem{corollary}[theorem]{Следствие}
  
  \newtheorem{definition}{Определение}[section]
  \newtheorem{utv}{Утверждение}[section]

\begin{document}
\newpage
\begin{minipage}{17cm}
    \begin{center}
        \textbf{\Large{Цена симметрии в связных графах}}
    \end{center}
\end{minipage}
\newline

\begin{minipage}{17cm}
    \begin{center}
        Михаил С. Терехов
    \end{center}
\end{minipage}
\newline

\begin{minipage}{15em}

\end{minipage}
$~~~~~~~~~~$
\begin{minipage}{13.1cm}
\footnotesize{В работе даётся ответ на вопрос, поставленный в совместной работе А. А. Клячко и Н. М. Луневой, об оптимальности оценки на цену симметрии в графах. Оригинальная оценка гласит, что, если в связном графе G можно удалить n вершин так, чтобы в нём не осталось связного подграфа изоморфного Г, то можно удалить не более n|V(Г)| вершин, образующих инвариантное относительно всех автоморфизмов графа G множество так, чтобы в графе не осталось подграфа изоморфного Г. Мы докажем, что существует граф Г, для которого эта оценка не является оптимальной.}

\end{minipage}

\section{Введение}
Рассмотрим следующую задачу.
\newline

\begin{minipage}{15em}

\end{minipage}
$~~~$
\begin{minipage}{15cm}
\normalsize{
	Пусть у нас есть некоторое множество шахматистов, и мы хотим провести турнир, в котором каждый игрок сыграет с каждым ровно по одному разу, но так,
чтобы не было подмножества из $k$ человек, попарно знающих друг друга (чтобы избежать большого количества договорных партий). Тогда если известно, что можно удалить $n$ игроков, чтобы условие выполнялось, то можно удалить не более $kn$ игроков, чтобы условие выполнялось, но <<честно>> 
--- то есть множество удаляемых игроков должно быть инвариантно относительно автоморфизмов на <<графе дружбы>>.
}
\end{minipage}
$~~~$
\newline
$~~~$
\newline
$~~~$
\newline
Более общий подход к данной задаче, заключён в следующей теореме (фактически именно это было доказано в [KlLu21], но для удобства читателя мы повторим доказательство в последнем разделе).
\newline
~~~~~
\newline
\textbf{Теорема 1.1 }[KlLu21]\textbf{.}
\textit{
Пусть группа $G$ действует на множестве $U$ и $\F~-~G$-инвариантное семейство
конечных подмножеств множества $U$, мощности которых ограничены в 
совокупности, а $X$ в $U$ --- конечная система представителей для этого 
семейства (то есть $X \cap F \neq \varnothing$, для всех $F \in \F$). Тогда найдётся 
$G$-инвариантная система представителей Y такая, что $|Y|  \leqslant |X|\max\limits_{F\in \F}|F|.$
При этом в качестве $Y$ можно взять следующее множество: 
$
Y=\left\{y\in U\;\Bigm|\;
|Gy\cap X| \geqslant {1\over \max\limits_{F\in\F}|F|}|Gy|\right\}
$.  
}

Слово \textit{семейство} здесь понимается как неупорядоченное
семейство, то есть $\F$ --- это просто некоторое множество подмножеств
множества $U$. Инвариантность семейства $\F$ следует понимать
естественным образом:
${gF\:=\{gf\;|\; f\in F\}\in\F}$ для всех $g\in G$ и $F\in\F$.

Из теоремы 1.1 несложно получается наш случай для графов.
\newline
~~~~~
\newline
\textbf{Следствие 1.1 }[KlLu21]\textbf{.}
\textit{
Пусть $\Gamma$ --- граф и $G$ --- конечный граф. Тогда, если в графе $\Gamma$ можно выбрать конечное множество вершин $X$ так,
чтобы каждый подграф графа $\Gamma$, изоморфный графу $G$, имел хоть одну
вершину из $X$, то
в графе $\Gamma$ можно выбрать конечное множество вершин $Y$, инвариантное
относительно всех автоморфизмов графа~$\Gamma$, так, чтобы опять
каждый подграф графа~$\Gamma$, изоморфный графу $G$, имел хоть одну
вершину из $Y$, причём
$|Y| \leqslant |X|\cdot(\hbox{\rm число вершин графа $G$})$.
}
\newline

Здесь слово \textit{граф} мы понимаем как неориентированный граф. Граф может иметь кратные рёбра и петли, или не иметь.

Как несложно убедиться, оценка в данном следствии может являться равентсвом. Например, если $\Gamma \simeq K_{|V(G)|}$, где $V(G)$ --- это множество вершин в графе $G$, а $K_n$ --- это полный граф на $n$ вершинах.

Несмотря на полученную оценку, логично задаться вопросом, нельзя ли её улучшить, скажем, на больших графах? Так мы приходим к определениям.

\textit{Вершинной представительностью} называют [KlLu21]  $\repV(K,\Gamma)$
графа $K$ в графе $\Gamma$ минимальное число $n$ такое, что в
графе~$\Gamma$ найдётся множество вершин $X$ мощности $n$,
удовлетворяющее следующему условию:
$$
\hbox{\sl каждый подграф графа $\Gamma$, изоморфный $K$, содержит вершину
из $X$.}
\eqno{(*)}
$$

\textit{Симметричной вершинной представительностью} называют [KlLu21] 
$\repVSym(K,\Gamma)$ графа $K$ в графе $\Gamma$ минимальное
число $n$ такое, что в графе $\Gamma$ найдётся инвариантное относительно
всех автоморфизмов множество вершин $X$ мощности $n$ с условием $(*)$.
\newline

Граф $K$ называют [KlLu21] 
\textit{вершинно дорогим}, если
$$
\hbox{
$\forall m \in Z$ найдётся граф
$\Gamma_m$ такой, что
}
\repVSym(K,\Gamma_m)=
\repV(K,\Gamma_m)\cdot(\hbox{|V($K$)|}) \geqslant m.
\eqno{(**)}
$$

Но даже так, на больших графах оценка будет неулучшаема, что легко понять, рассмотрев дизъюнктное объединение $m$ полных графов $K_{|V(G)|}$.

Однако иногда логично наложить ограничение на исходный граф. Например, в задаче про шахматистов можно предположить, что мир шахмат тесен, и наш <<граф дружбы>>$~$связен. 

Вершинно дорогой граф $K$ называют [KlLu21]\textit{\(вершинно\) дорогим в классе
графов $\cal K$}, если графы $\Gamma_m$ в $(**)$
могут выбраны из класса~$\cal K$.

Таким образом, мы приходим к основному вопросу, которому посвящена данная работа.
\newline
~~~~~
\newline
\textbf{Вопрос 1} [KlLu21]\textbf{.}
Верно ли, что любой конечный связный граф является дорогим в классе
связных графов?
\newline

В работе [KlLu21] доказано, что любой связный граф без висячих вершин является вершинно дорогим в классе связных графов.
Также доказано, что все связные графы с менее чем пятью вершинами являются вершинно дорогими в классе связных графов.

Так мы получаем, что самый маленький по количеству вершин претендент на роль вершинно недорогого связного графа это $D_5$.

$~~~~~~~~~~~~~~~~~~~~~~~~~~~~~~~~~~~~~~~~~~~~~~~$
\unitlength 1.00mm
\linethickness{0.4pt}
\begin{picture}(36.00,29.00)
\put(5.00,8.00){\line(1,1){10.00}}
\put(15.00,18.00){\line(1,0){20.00}}
\put(15.00,18.00){\circle*{2.00}}
\put(25.00,18.00){\circle*{2.00}}
\put(35.00,18.00){\circle*{2.00}}
\put(5.00,8.00){\circle*{2.00}}
\put(5.00,28.00){\line(1,-1){10.00}}
\put(5.00,28.00){\circle*{2.00}}
\put(20.00,4.00){\makebox(0,0)[cc]{$D_5$}}
\put(20.00,-1.00){\makebox(0,0)[cc]{Рис. 1}}

\end{picture}
\newline

Для этого графа была доказана теорема 1.2.
\newline
~~~~~
\newline
\textbf{Теорема 1.2 }[KlLu21]\textbf{.}
\textit{
Граф $D_5$ не является дорогим в классе вершинно
транзитивных связных графов.
Более точно,
если $\Gamma\supseteq D_5$ --- вершинно транзитивный неориентированный
связный граф,
имеющий больше пяти вершин, и представительность 
$\repV(D_5,\Gamma)$ конечна, то $\repVSym(D_5,\Gamma) < 5 \repV(D_5,\Gamma)$.
}
\newline

Мы ответим на Вопрос 1 отрицательно, доказав  более сильный факт (см. Теорему 2.2), из которого 
легко будет следовать следующая теорема. 
\newline
~~~~~
\newline
\textbf{Теорема 1.3.}
Граф $D_5$ не является вершинно дорогим в классе связных графов.
\newline
$~~~~~~~~~~~~~~~~~~~~~~~~~~~~~~~~~~~~~~~~~~~~~~~$
\newline
Выражаю благодарность А.А. Клячко и А.Л. Таламбуце за ценные замечания. 
\section{Доказательство теоремы 1.3}



\textit{Орбитой вершины v} графа $\Gamma = (V,E)$ назовём множество вершин $V_1(v) = 
\left\{g(v)\Bigm|~g\in Aut(\Gamma) \right\}$.
\newline
Очевидно, что множество вершин графа распадается на орбиты: $V(\Gamma) = V_1 \sqcup...\sqcup V_k$.


В каждом графе $\Gamma$ выберем множество вершин $M(\Gamma)$ так, чтобы каждый подграф, изоморфный графу 
$D_5$, содержащийся в $ \Gamma$, имел хотя бы одну вершину из множества $M(\Gamma)$ при этом $|M(\Gamma)|=\repV(D_5,\Gamma).$ Скажем, что множество $M(\Gamma)$ \textit{отмечено} в графе
$\Gamma$.
В каждом графе $\Gamma$ выберем множество вершин $M_{sym}(\Gamma)$ так, чтобы каждый подграф, изоморфный графу 
$D_5$, содержащийся в $ \Gamma$, имел хотя бы одну вершину из множества $M_{sym}(\Gamma)$ при этом $|M_{sym}(\Gamma)|=\repVSym(D_5,\Gamma)$, а также $M_{sym}(\Gamma)$ являлось инвариантным множеством, относительно $Aut(\Gamma)$. Скажем, что множество  $M_{sym}(\Gamma)$ 
\textit{симметрично отмечено} в графе $\Gamma$. Будем считать, что $M_{sym}(\Gamma)$ выбирается, исходя из применения теоремы 1.1 к множеству $M(\Gamma)$.

\utv
Пусть связный граф $\Gamma$ $\not\simeq$ $K_5$ такой, что $5\repV(D_5,\Gamma)=\repVSym(D_5,\Gamma)>0$. Тогда 
в графе есть хотя бы две орбиты.
\proof
Для одной орбиты мы получаем вершинно транзитивный связный граф, содержащий подграф, изоморфный графу $D_5$. В этом случае по теореме 1.2 выполнено неравенство $\repVSym(D_5,\Gamma) < 5\repV(D_5,\Gamma)$.
\endproof
Будем говорить, что орбита $A$ \textit{пересекается} с $D_5$, если граф $\Gamma$ содержит подграф, изоморфный графу 
$D_5$, у которого хотя бы одна вершина лежит в $A$. 
\newline

\begin{lemma}

Пусть связный граф $\Gamma$ $\not\simeq$ $K_5$ такой, что $5\repV(D_5,\Gamma)=\repVSym(D_5,\Gamma)>0$. Тогда для каждой орбиты $A$ верно, что
$|A \cap M(\Gamma)| \in \{0,\frac{1}{5}|A|\}.$ При этом, если орбита не пересекает
$D_5$, то $|A \cap M(\Gamma)| = 0$, а если орбита $A$ 
пересекает $D_5$, то 
$|A \cap M(\Gamma)|=\frac{1}{5}|A|$.

\proof
Из теоремы 1.1, применённой к семейству $\F = \{$ все подграфы в  $\Gamma$, изоморфные графу $D_5 \},~X=M(\Gamma)$ и $G=Aut(\Gamma)$, сразу же вытекает, что для любой орбиты A графа $\Gamma$ либо $|A\cap M(\Gamma)|=\frac{1}{5}|A|$ и $A \subseteq M_{sym}(\Gamma)$
, либо $A\cap M(\Gamma) = \varnothing = A\cap M_{sym}(\Gamma)$.
\newline
Теперь остаётся понять, почему орбита $A$, пересекающая $D_5$, содержит отмеченную вершину.
Если орбита содержит подграф, изоморфный графу $D_5$, то доказывать нечего. Допустим, что нет.
Тогда применяя Теорему 1.1 к 
\newline
$$U = V(\Gamma) \setminus A,~\F = \{V(H)\setminus A~|~H \simeq D_5,~H \cap A \neq \varnothing \},~X=M(\Gamma)~\text{и}~ G=Aut(\Gamma),$$
\newline 
мы получаем, что симметричное множество представителей $Y$ пусто (поскольку в каждой орбите меньше четверти вершин лежит в $X$), 
что является очевидным противоречием.
\endproof
\end{lemma}

\begin{lemma}
Пусть между некоторыми орбитами $A$ и $B$ есть хотя бы одно ребро, и 
$ S_1 \subseteq A$ --- некоторое подмножество вершин. Тогда для множества вершин $ S_2 \subseteq B$, связанных хотя бы одним ребром с множеством $S_1$, выполнено неравенство $|S_2| \geqslant |S_1|\frac{|B|}{|A|}$.
\proof
Обозначим $\frac{|S_1|}{|A|}$ за $c$.
Пусть из каждой вершины в $A$, рёбер, ведущих в орбиту $B$, ровно $k$ штук. Тогда из каждой вершины в $B$ ведёт $k \frac{|A|}{|B|}$ рёбер в орбиту $A$. 
Теперь предположим, что все рёбра ведущие из множества $S_1$ в орбиту $B$ ведут в подмножество $S_2$. Тогда рёбер, ведущих из $S_1$ в $S_2$, ровно $ck|A|$, значит степень каждой вершины в $S_2$ хотя бы $\frac{c k|A|}{|S_2|}$, но степень каждой вершины из $S_2$ равна  
$k \frac{|A|}{|B|}$, значит $\frac{c k|A|}{|S_2|} \leqslant k \frac{|A|}{|B|}$, следовательно 
$|S_2| \geqslant c|B|=|S_1|\frac{|B|}{|A|}$.
\endproof
\end{lemma}

\begin{lemma} 
Пусть связный граф $~\Gamma$ $\not\simeq$ $K_5$ такой, что  
$5\repV(D_5,\Gamma)=\repVSym(D_5,\Gamma)>0$. Тогда никакая орбита графа $\Gamma$ не содержит подграф, изоморфный графу $D_5$.
\newline
\proof
Сначала заметим, что компоненты связности в одной орбите между собой изоморфны. 
В самом деле, переведя автоморфизмом вершину из одной связной компоненты в другую компоненту, связная компонента перейдёт в связную комоненту и установит изоморфизм между компонентами.
Теперь будем рассуждать от противного. Предположим, что некоторая орбита $A$ графа$~\Gamma$ содержит подграф, изоморфный графу $D_5$. Рассмотрим 2 случая:
\newline
1) Cвязная компонента орбиты $A$ не изоморфна $K_5$.

Тогда по теореме 1.2 для каждой связной компоненты $A_i$ орбиты $A$ выполнено неравенство 
\newline
$|M(\Gamma)\cap V(A_i)|>\frac{1}{5}|V(A_i)|$, следовательно 
$|M(\Gamma)\cap V(A)|>\frac{1}{5}|V(A)|$, но по лемме 2.1 в орбите $A$ отмечено не более $\frac{1}{5}|V(A)|$ вершин.
\newline
2) Cвязная компонента орбиты $A$ изоморфна $K_5$.

Тогда чтобы каждый подграф, изоморфный графу $D_5$, имел вершину из 
$M(\Gamma)$, в каждой компоненте изоморфной $K_5$ должна быть хотя бы одна вершина из $M(\Gamma)$. А значит, ровно одна, в силу  
леммы 2.1. 
Теперь рассмотрим орбиту $B$ из которой есть ребро в орбиту $A$.
По лемме 2.2 хотя бы $\frac{4}{5}|B|$ вершин 
имеют общее ребро хотя бы с одной из $\frac{4}{5}|A|$ неотмеченных вершин в $A$. Каждая вершина из $B$, связанная с неотмеченной вершиной из орбиты $A$ принадлежит подграфу, изоморфному графу $D_5$, в котором все остальные вершины неотмечены (рис. 2). 
Значит, все вершины из $B$, которые имеют общее ребро с неотмеченной вершиной из $A$ должны быть отмечены в $\Gamma$. Но по лемме 2.1 в орбите $B$ отмечено не более $\frac{1}{5}|B|$ вершин.
\newline
\unitlength 1.00mm
\linethickness{0.4pt}
\begin{picture}(116.00,50.00)

\put(50.00,10.00){\circle{2.00}}
\put(70.00,10.00){\circle{2.00}}
\put(60.00,40.75){\circle*{2.00}}
\put(43.80,29.00){\circle{2.00}}
\put(76.20,29.00){\circle{2.00}}
\put(97.00,29.00){\circle*{2.00}}

\put(45.00,29.00){\line(1,0){30.0}}
\put(77.00,28.90){\line(1,0){19.0}}
\put(77.00,29.00){\line(1,0){19.0}}
\put(77.00,29.20){\line(1,0){19.0}}
\put(77.00,29.10){\line(1,0){19.0}}

\put(44.70,29.00){\line(5,4){14.2}}

\put(44.80,28.60){\line(4,-3){24.0}}
\put(45.10,28.60){\line(4,-3){24.0}}
\put(44.60,28.60){\line(4,-3){24.0}}
\put(45.0,28.60){\line(4,-3){24.0}}
\put(44.90,28.60){\line(4,-3){24.0}}

\put(51.00,10.00){\line(1,0){18.0}}
\put(51.00,10.10){\line(1,0){18.0}}
\put(51.00,10.20){\line(1,0){18.0}}
\put(51.00,10.30){\line(1,0){18.0}}
\put(51.00,10.40){\line(1,0){18.0}}

\put(49.58,10.93){\line(1,3){9.67}}
\put(70.42,10.93){\line(-1,3){9.67}}
\put(75.20,29.00){\line(-5,4){14.2}}

\put(49.50,10.9){\line(-1,3){5.7}}

\put(70.10,11.00){\line(1,3){5.7}}
\put(69.90,11.00){\line(1,3){5.7}}
\put(70.00,11.00){\line(1,3){5.7}}
\put(70.20,11.00){\line(1,3){5.7}}

\put(97.00,32.00){\makebox(0,0)[cc]{$b$}}
\put(46.00,10.00){\makebox(0,0)[cc]{$a_4$}}
\put(56.00,42.75){\makebox(0,0)[cc]{$a_2$}}
\put(40.00,30.00){\makebox(0,0)[cc]{$a_3$}}
\put(76.00,32.00){\makebox(0,0)[cc]{$a_1$}}
\put(74.00,10.00){\makebox(0,0)[cc]{$a_5$}}

 \put(84.00,4.00){\makebox(0,0)[cc]{Рис. 2}}

\end{picture}
\newline
На рисунке вершины $a_1, a_2, a_3, a_4, a_5$ принадлежат орбите $A$ 
при этом $a_2\in M(\Gamma)$. 
\newline
Вершина $b$ принадлежит  
орбите $B$, и пересекается с $D_5$, а значит является отмеченной.
\end{lemma}
\begin{corollary}
Пусть связный граф $\Gamma$ $\not\simeq$ $K_5$ такой, что $5\repV(D_5,\Gamma)=\repVSym(D_5,\Gamma)>0$. Тогда $\repVSym(D_5,\Gamma)<$ $|V(\Gamma)|$. 
\proof Пусть $A~-$ некоторая орбита. Тогда по лемме 2.3, 
любой подграф, изоморфный графу $D_5$, имеет хотя бы вершину из
$\Gamma \setminus A$, следовательно $\repVSym(D_5,\Gamma)\leqslant |V(\Gamma\setminus A)|<|V(\Gamma)|$.
\endproof

\end{corollary}
\begin{definition}
Орбиту $A$ назовём отмеченной, если $A\subseteq M_{sym}(\Gamma)$ и неотмеченной если $A \cap M_{sym}(\Gamma) = \varnothing$.
\end{definition}
Любая орбита либо отмеченная, либо неотмеченная.

\begin{lemma}
Пусть связный граф $\Gamma$ $\not\simeq$ $K_5$ такой, что $5\repV(D_5,\Gamma)=\repVSym(D_5,\Gamma)>0$. Тогда $M(\Gamma) \subseteq M_{sym}(\Gamma)$.
\proof
Любая отмеченная орбита $A$ пересекается с $D_5$(иначе её можно было бы не отмечать), значит по лемме 2.1 
\newline
$|A \cap M(\Gamma)|=\frac{1}{5}|A|$, следовательно, в силу равенства 
$5\repV(D_5,\Gamma)=\repVSym(D_5,\Gamma)$, в неотмеченных орбитах мы не можем отмечать вершины.
\endproof
\end{lemma}
\begin{lemma}
Для связных графов $\Gamma_1,\Gamma_2,...,\Gamma_k$
выполнено равенство 
$|M_{sym}(\Gamma_1 \sqcup \Gamma_2 \sqcup ...\sqcup \Gamma_k)| =$
\newline
$=|M_{sym}(\Gamma_1)|+|M_{sym}(\Gamma_2)|+...+|M_{sym}(\Gamma_k)|$.
\proof
Нетрудно понять, что достаточно доказать утверждение для случая, когда 
$\Gamma_1 \simeq \Gamma_2 \simeq ...\simeq \Gamma_k $.
Пусть $V(\Gamma_1)=V_{11}\sqcup\ldots \sqcup V_{1s}$, где
$V_{1r}~-$ орбиты графа $\Gamma_1$. Тогда в силу 
$\Gamma_1 \simeq \Gamma_i$, можно считать что 
$V(\Gamma_i)=V_{i1}\sqcup\ldots \sqcup V_{is}$ при этом для 
любого изоморфизма $\alpha:\Gamma_1\rightarrow\Gamma_i$ 
верно что $\alpha(V_{1r})=V_{ir}$. Предположим теперь, 
что $M_{sym}(\Gamma_1)=V_{11}\sqcup\ldots V_{1j}$. 
Тогда в качестве $M_{sym}(\Gamma_1 \sqcup \Gamma_2 \sqcup ...\sqcup \Gamma_k)$ подойдёт множество 
$V_{11}\sqcup\ldots V_{1j} \sqcup\ldots\sqcup V_{k1}\sqcup\ldots V_{kj}$.
\endproof
\end{lemma}

\theorem
Для любого связного графа $\Gamma$ $\not\simeq$ $K_5$, содержащего
подграф, изоморфный графу $D_5$, выполнено неравенство: $5\repV(D_5,\Gamma)$ $>\repVSym(D_5,\Gamma)$.
\proof

Предположим, утверждение теоремы 2.2 неверно, тогда рассмотрим контрпример, минимальный по количеству вершин.
По следствию 2.1 в $\Gamma$ есть неотмеченная орбита.
\newline 
Рассмотрим неотмеченную орбиту $A$, которая соединена ребром
 с отмеченной орбитой $B$.
По лемме 2.4 имеем $|A\cap M(\Gamma)|=0$, значит по лемме 2.1 получаем, что орбита $A$ не пересекается с $D_5$. Теперь выбросим орбиту $A$. 
Далее докажем цепочку неравенств 
$$ \repVSym(D_5,\Gamma) \stackrel{(1)}{\leqslant} \repVSym(D_5,\Gamma\backslash A)  \stackrel{(2)}{\leqslant} 5\repV(D_5,\Gamma\backslash A)  \stackrel{(3)}{=}5\repV(D_5,\Gamma)\stackrel{(4)}{=}\repVSym(D_5,\Gamma).$$

Неравенство $(1)$ выполнено в силу того, что 
инвариантное множество вершин относительно 
 $Aut(\Gamma\backslash A)$ является инвариантным относительно $Aut(\Gamma)$ и того, что вершины из орбиты $A$ не отмечены. 

Неравенство $(2)$ выполнено в силу следствия 1.1. 

Равенство $(3)$ выполнено в силу того, что
$D_5$ не пересекает $A$.

Равенство $(4)$ выполнено в силу того, что
$\Gamma$ -- гипотетический контрпример к теореме 2.2.

Так мы получаем, что неравенство (2) является равенством. 
\newline
Обозначим компоненты связности $\Gamma\backslash A$ как 
$\Gamma_1,\Gamma_2,...,\Gamma_k$, тогда по лемме 2.5 и следствию 1.1. выполнено 
\newline
 $$\repVSym(D_5,\Gamma\backslash A) =
 \repVSym(D_5,\Gamma_1)+...+\repVSym(D_5,\Gamma_k) \leqslant 5\repV(D_5,\Gamma_1)+...+5\repV(D_5,\Gamma_k)= 5\repV(D_5,\Gamma\backslash A).$$ Значит, для всех $\Gamma_i$ выполнено равенство 
$\repVSym(D_5,\Gamma_i)=5\repV(D_5,\Gamma_i)$, откуда в силу минимальности контрпримера следует, что либо $\Gamma_i$ не содержит подграфа изоморфного $D_5$, либо $\Gamma_i \simeq K_5$.
Но так как ни один $D_5$ не пострадал при удалении орбиты $A$, то вершины из орбиты $B$ по-прежнему будут пересекаться с $D_5$, значит вершины из орбиты $B$ будут принадлежать компоненте изоморфной $K_5$.  
В этом случае орбита $A$ пересекается с $D_5$ (рис. 3), следовательно по лемме 2.1, $|A \cap M(\Gamma)|=\frac{1}{5}|A|$, но в то же время по лемме 2.4 выполнено
$M(\Gamma) \subseteq M_{sym}(\Gamma)$ и по нашему предположению 
орбита $A$ является неотмеченной, то есть $A \cap M_{sym}(\Gamma) = \varnothing$.
\endproof

\unitlength 1.00mm
\linethickness{0.4pt}
\begin{picture}(116.00,50.00)

\put(50.00,10.00){\circle{2.00}}
\put(70.00,10.00){\circle{2.00}}
\put(60.00,40.75){\circle{2.00}}
\put(43.80,29.00){\circle{2.00}}
\put(76.20,29.00){\circle{2.00}}
\put(97.00,29.00){\circle{2.00}}

\put(45.00,29.00){\line(1,0){30.0}}
\put(77.00,28.90){\line(1,0){19.0}}
\put(77.00,29.00){\line(1,0){19.0}}
\put(77.00,29.20){\line(1,0){19.0}}
\put(77.00,29.10){\line(1,0){19.0}}

\put(44.70,29.00){\line(5,4){14.2}}

\put(44.80,28.60){\line(4,-3){24.0}}
\put(45.10,28.60){\line(4,-3){24.0}}
\put(44.60,28.60){\line(4,-3){24.0}}
\put(45.0,28.60){\line(4,-3){24.0}}
\put(44.90,28.60){\line(4,-3){24.0}}

\put(51.00,10.00){\line(1,0){18.0}}
\put(51.00,10.10){\line(1,0){18.0}}
\put(51.00,10.20){\line(1,0){18.0}}
\put(51.00,10.30){\line(1,0){18.0}}
\put(51.00,10.40){\line(1,0){18.0}}

\put(49.58,10.93){\line(1,3){9.67}}
\put(70.42,10.93){\line(-1,3){9.67}}
\put(75.20,29.00){\line(-5,4){14.2}}

\put(49.50,10.9){\line(-1,3){5.7}}

\put(70.10,11.00){\line(1,3){5.7}}
\put(69.90,11.00){\line(1,3){5.7}}
\put(70.00,11.00){\line(1,3){5.7}}
\put(70.20,11.00){\line(1,3){5.7}}
\put(101.00,32.00){\makebox(0,0)[cc]{$a\in A$}}
\put(80.00,32.00){\makebox(0,0)[cc]{$b\in B$}}
\put(76.00,4.00){\makebox(0,0)[cc]{Рис. 3}}

\end{picture}
\newline
На рисунке жирным отмечены рёбра $D_5$, пересекающего орбиту $A$.
\newline
Из теоремы 2.2 легко вытекает теорема 1.3.
\newline

\section{Доказательство теоремы 1.1}
\textit{Доказательство} {[KlLu21].}
Положим $m=\max\limits_{F\in\F}|F|$ и рассмотрим следующее множество

$
Y=\left\{y\in U\;\Bigm|\;
|Gy\cap X| \geqslant {1\over m}|Gy|\right\}
$
(в частности, $Y$ не содержит точек с бесконечной орбитой).
Ясно, что это множество $G$-инвариантно. Ясно также, что
$|Y| \leqslant m|X|$ (поскольку для каждой орбиты $Gu$ имеет место
неравенство $|Gu\cap Y| \leqslant m|Gu\cap X|$).

Осталось показать, что $Y$ является системой представителей для $\F$.
Возьмём какое-то множество $F\in\F$.
Каждое множество $gF$ (где $g\in G$) принадлежит $\F$ в силу
инвариантности семейства $\F$ и, следовательно, пересекается с $X$.
Значит,
$$
G=\bigcup_{f\in F}\{g\in G\;|\;gf\in X\}.
$$
Каждое из множеств $\{g\in G\;|\;gf\in X\}$ является
либо пустым, либо объединением конечного числа левых смежных
классов группы $G$ по стабилизатору $St(f)$ точки $f$:
$$
\{g\in G\;|\;gf\in X\}=
\bigcup_{x\in X}\{g\in G\;|\;gf=x\}=$$

$$=\bigcup_{x\in X\cap Gf}g_x\cdot St(f),~
\{\text{где } g_x\in G \text{ фиксированы так, что } g_xf=x\}.
$$
Таким образом, мы получили разложение группы~$G$ в конечное объединение
левых смежных классов по некоторым подгруппам.
Воспользуемся теперь теоремой Б.~Неймана [Neu54] (утверждение~4.5):
если группа $G$ покрывается конечным числом смежных классов по некоторым
\(необязательно разным\) подгруппам: $G=g_1G_1\cup\dots\cup g_sG_s$,
то
\
$\displaystyle\sum {1\over|G:G_i|} \geqslant 1$
\(где обратный к бесконечному кардиналу
считается нулём\).
Следовательно, (учитывая то, что индекс стабилизатора равен длине орбиты)
мы получаем
$$
1 \leqslant \sum_{f\in F}{1\over|G: St(f)|}\cdot|Gf\cap X|=
\sum_{f\in F}{|Gf\cap X|\over|Gf|}.
$$
Поскольку число слагаемых в этой сумме равно $|F|\leqslant m$, по крайней мере
одно из слагаемых должно быть не меньше чем $1/m$, то есть
${|Gf\cap X|/|Gf|} \geqslant 1/m$, что означает $f\in Y$ (по определению
множества $Y$) и завершает доказательство.
\endproof

СПИСОК ЦИТИРОВАННОЙ ЛИТЕРАТУРЫ
\newline
[KlLu21]	~~~A. A. Klyachko, N. M. Luneva, Invariant systems of representatives, or The cost of symmetry,
Discrete Mathematics, 344:6 (2021), 112361. См. также arXiv:1908.03315.
\newline
[Neu54]		~~~~~~B.H. Neumann, Groups covered by 
permutable subsets, J. London Math. Soc., s1-29:2 (1954), 236-248.

\end{document}